\documentclass[11pt,twoside]{article}
\usepackage{amssymb}
\usepackage[mathscr]{eucal}
\usepackage{eufrak}
\usepackage{amsmath,amsthm}
\usepackage{mathrsfs}
\usepackage[english]{babel}
\usepackage{lgreek}
\usepackage[shortlabels]{enumitem}
\usepackage[colorlinks=true,
   urlcolor=blue,           filecolor=green,      
   citecolor=green,      
   linkcolor=red,           bookmarks=true,
  unicode,
   plainpages=false,   ]{hyperref}
\usepackage{color}

\definecolor{royalblue}{rgb}{0,0,0.128}
\def\bp{\begin{proof}}
\def\ep{\end{proof}}
\def\n{\nabla}

\def\sfrac#1#2{\mbox{\Large$\frac{#1}{#2}$}}

\def\intll#1#2{\int\limits_{#1}^{#2}}
\def\dm{|\hskip-0.05cm|}

\def\OO{\Omega}
\def\displ{\displaystyle}
\def\VSE{\vspace{6pt}\\&\displ }
\def\VS{\vspace{6pt}\\\displ }
\def\rf#1{{\rm(\ref{#1})}}

\def\R{\Bbb R}
\def\N{\Bbb N}
\def\A{\Bbb A}
\def\à{à}

\def\be{\begin{equation}}
\def\ba{\begin{array}}
\def\ea{\end{array}}
\def\ee{\end{equation}}
\def\vs1{\vspace{1ex}}

\def\po{\partial\Omega}

\def\�{\'{e}}
\def\�{\`{e}}

\newtheorem{lemma}
{\bf Lemma} 
\font\sc=cmcsc10
\title{\Large\bf A new proof of the \textbf{Th\'eor\�me de Structure} related to\\a weak solution to the Navier-Stokes equations}
\author{{\it In memoria del Professor Herman Sohr}\vspace{6pt}\\\sc  Paolo Maremonti and Filippo Palma
\thanks{Dipartimento di Matematica e Fisica,  
Universit\`{a} degli 
Studi della Campania
``L. Vanvitelli'', via Vivaldi 43, 81100 \null\hskip0.55cmCaserta,
 Italy.\newline\null\hskip0.55cm
paolo.maremonti@unicampania.it \;\;and\;\; filippo.palma@unicampania.it\newline\null\hskip0.55cm The  research activity  is performed under the
auspices of   GNFM-INdAM. }}
\date{}
\begin{document}
\markboth{\footnotesize\rm    P. Maremonti and F. Palma} {\footnotesize\rm A new proof of the {\it  Th\'{e}or\`{e}me de Structure} related to a weak solution to the Navier-Stokes equations
}
\maketitle 
{\small\bf Abstract} - {\small It is well known that a Leray's weak solution to the Navier-Stokes Cauchy problem enjoys a partial regularity which is known in the literature as the  {\it Th\'eor\�me de Structure} of
 a Leray's  weak solution. As well, this result has been extended by some authors  to the case of the IBVP. \newline\null\hskip0.5cm In this note,    we achieve   the {\it Th\'eor\�me de Structure}  by means of a new  proof. Our proof is based on {\it a priori} estimates for a suitable approximating sequence. In this way our result covers a more general setting in the sense that,  {\it e.g.}, we can also include  the case of the weak solutions furnished by Hopf for an IBVP in bounded domains without requiring an energy inequality in a strong form,  but just employing a priori estimates on the Galerkin approximation.
}
\vskip0.2cm \par\noindent{\small Keywords: Navier-Stokes equations,   weak solutions, partial regularity. }
  \par\noindent{\small  
  AMS Subject Classifications: 35Q30, 35B65, 76D05.}  
\noindent
\newcommand{\red}{\protect\bf}
\renewcommand\refname{\centerline
{\red {\normalsize \bf References}}}
\newtheorem{ass}
{\bf Assumption} 
\newtheorem{defi}
{\bf Definition} 
\newtheorem{tho}
{\bf Theorem} 
\newtheorem{rem}
{\sc Remark} 
\newtheorem{coro}
{\bf Corollary} 
\newtheorem{prop}
{\bf Proposition} 
\renewcommand{\theequation}{\arabic{equation}}
\setcounter{section}{0}
\section{Introduction}
In this note we consider the initial boundary value problem (IBVP) for the Navier-Stokes equations:
\be\label{NS}\ba{cc}v_t+v\cdot\n v+\n\pi=\Delta v\,,&\n\cdot v=0\,,\mbox{ in }(0,T)\times \OO\,,\VS v=0\mbox{ on }(0,T)\times\po\,,& v=v_0\mbox{ on }\{0\}\times\OO\,,\ea\ee where $\OO\subseteq\R^3$ can be a   bounded or an exterior  domain, whose boundaries $\po$ for simplicity 
are assumed smooth, a half-space or the whole space. We set $w_t:=\frac{\partial w} {\partial t}$ and $w\cdot\n w:= (w\cdot\n) w$. \par In \cite{L}, J. Leray established the first result related to a weak solution to the Navier-Stokes Cauchy problem. The peculiarity of the result is the fact that the existence theorem  furnishes a solution global in time (for all $t>0$) with no restriction on the size of  the initial datum that for simplicity belongs to $L^2(\R^3)$ with divergence free in a weak sense.\par  Contextually, he implicitly posed the question of the regularity of such a solution. For this question, Leray  proposed the so called {\it Th\'eor\�me de Structure of a weak solution} (in that time for his weak solution). \par{\it The structure theorem states that the weak solution $v$ is regular at least on a sequence of open intervals $(\theta_\ell,T_\ell)$, $\ell\in\N$, and one of these intervals is of the kind $(\theta,\infty)$ with $\theta\leq c\dm v_0\dm_2^4$, and the Lebesgue's measure of $(0,\theta)-{\underset{\ell\in\N}\cup}(\theta_\ell,t_\ell)$ is null.}\par In the following  $v$ denotes  a weak solution to problem \rf{NS}. \par Following J. Leray, the proof of the {\it Th\'eor\�me de Structure} is based on the following {\it``ingredients''}. \begin{itemize}\item[i.] The validity of the energy inequality in strong form, that is,
\be \label{ES}\dm v(t)\dm_2^2+2\intll st\dm\n v(\tau)\dm_2^2d\tau\leq \dm v(s)\dm_2^2,\,\mbox{ a.e. in }{ s\in (0,t) \mbox{ and for }s=0.}\ee Set ${\mathcal T}_L:=\{s:\mbox{ energy inequality \rf{ES} holds}\}$, one gets\be\label{PLIM}\mbox{for all }s\in{\mathcal T}_L\,,\;\lim_{t\to s^+}\dm v(t)-v(s)\dm_2=0\,.\ee  \item[ii.]  The possibility that for $v(s,x)\in J^{1,2}(\OO)$ there exists a regular solution $w(t,x)$ on a maximal interval $(s,s+T(s))$, and there exists an $\eta_0>0$ such that if $\dm v(\theta)\dm_2\dm\n v(\theta)\dm<\eta_0$,  then the solution $w(t,x)$ is regular for all $t\in(\theta,\infty)$.\item[iii.] The uniqueness between the weak solution $v(t,x)$ and the regular solution $w(t,x)$ corresponding to $v(s,x)\in J^{1,2}(\OO)$ holds on $(s,s+T(s))$, or $(\theta,\infty)$, provided that $s\in {\mathcal T}_L$.\end{itemize}  
\par In \cite{SW} another approach is proposed in order to get  the partial regularity in the sense of the {\it  Th\'eor\�me de Structure}. The idea is to break down a weak solution in the sum of two functions. Each function is the solution to a suitable linear problem. One is the  solution to the linear IBVP of the Stokes system with zero body force. Another is a solution to the linear IBVP of the Stokes system with homogeneous initial datum and body force equal to the convective term of the weak solution. The authors are able to prove a characterization of the validity of Leray's {\it Th\'eor\�me de Structure} in terms of the regularity of a weak solution.
\par In order to better state our chief result, we introduce some notations and the definition of the regular solution and  of a weak solution. \par {We denote by $$\mathscr C_0(\OO):=\{\varphi:\varphi\in C_0^\infty(\OO)\mbox{ and }\n\cdot\varphi=0\}\,.$$ We indicate by $J^2(\OO)$ and by $J^{1,2}(\OO)$ the completion of $\mathscr C_0(\OO)$  in $L^2(\OO)$ and in $W^{1,2}(\OO)$, respectively. For $q\in[1,\infty)$, by the symbol $L^q(0,T;X)$ we denote the classical Bochner space of measurable functions $u:t\in(0,T)\to X$, $X$-Banach space, and  with $\intll0T\dm u(t)\dm_X^qdt<\infty$.}
\begin{defi}\label{RS}{\sl For all $v_0\in J^{1,2}(\OO)$, the pair $(v,\pi)$ with $v:(0,T)\times\OO\to\R^3$ and $\pi:(0,T)\times\OO\to\R$ is said a regular solution to the IBVP \rf{NS} provided that $(v,\pi)$ satisfies the problem \rf{NS} almost everywhere in $(t,x)$,  and \be\label{RS}v\in C(0,T;J^{1,2}(\OO))\cap L^2(0,T;W^{2,2}(\OO))\mbox{\; and\; }v_t\,,\;\n\pi\in L^2(0,T;L^2(\OO))\,,\ee with $\displ\lim_{t\to0}\dm v(t)-v_0\dm_2=0$\,.  }\end{defi}
\begin{defi}\label{WS}{\sl For all $v_0\in J^2(\OO)$, a field $v:(0,\infty)\times\OO\to\R^3$ is said a weak solution corresponding to the datum $v_0$ if \begin{itemize}\item[{\rm1)}]for all $T>0$, $v\in L^\infty(0,T;J^2(\OO))\cap L^2(0,T;J^{1,2}(\OO))$,\item[{\rm 2)}]for all $T>0$ and  $t,s\in(0,T)$,  the field $v$
satisfies the integral   equation:
\newline \centerline{$\displ\intll
st\Big[(v,\varphi_\tau)-(\nabla
v,\nabla
\varphi)+(v\cdot\nabla\varphi,v)\Big]d\tau+(v(s),\varphi
(s))=(v(t),\varphi(t))$,} for all $\varphi(t,x)\in C_0^1([0,T)\times\OO)$ with $\n\cdot\varphi(t,x)=0$\,,\item[{\rm 3)}]$\displ\lim_{t\to0}\,(v(t),\varphi)=(v_0,\varphi)\,,\mbox{ for all }\varphi\in\mathscr C_0(\OO)\,.$\end{itemize}}\end{defi} \par The aim of this note is to state a different technique to achieve the   {\it Th\'eor\�me de Structure} of a weak solution. The major difference consists  in the fact that we do not   work on the weak solution $v$, but just considering a priori estimates achieved on the approximating sequence $\{v^m\}$ of the weak solution $v$. Hence, we employ no property of the weak solution  stated in the previous items i.-iii or like in \cite{SW}.
\par We are in a position to state our chief result on the structure theorem related to a weak solution to the initial boundary value problem \rf{NS}. 
\begin{tho}\label{CT}{\sl Let $\{v_0^m\}\subset\mathscr C_0(\OO)$ and $\dm v_0^m-v_0\dm_2\to0$.  Let $\{v^m\}$
be a sequence such that  \be\label{CT-I}\{v^m\}\subset C(0,T;J^{1,2}(\OO))\cap L^2(0,T;W^{2,2}(\OO)) \mbox{\; and\; }\{v^m_t\}\subset L^2(0,T;J^2(\OO)), \mbox{ for all }T>0.\ee   Assume that    for all $m\in\N$,  $v^m$ satisfies the energy relation \be\label{EW}\dm v^m(t)\dm_2^2+2\intll0t\dm\n v^m(\tau)\dm_2^2d\tau\leq \dm v_0^m\dm_2^2\,,\mbox{ for all }t>0\,,\ee and,  in the case of $\Omega$ unbounded, the following limit property uniformly in $m$ holds:
\be\label{CT-II}\lim_{R\to\infty}\intll0 t\dm v^m(\tau)\dm_{L^2(|x|>R)}d\tau=0\,,\mbox{ for all } t>0\,.\ee {Assume that, there exists\,\footnote{\,We assume $\alpha\geq\frac12$ for the sake of the brevity. Actually, it is enough to assume $\alpha>0$.} an $\alpha\geq\frac12$ such that, for all $T>0$,}    
\be\label{DS}{\intll0T\dm D^2v^m(t)\dm_2^\alpha dt\leq A(T,\dm v_0\dm_2)<\infty}\,,\mbox{ for all }m\in\N\,.\ee Assume that for all $m\in \N$ the following differential inequality for the Dirichlet norm of $v^m$ holds:
\be\label{DDN}\sfrac d{dt}\dm\n v^m(t)\dm_2^2+\dm P\Delta v^m(t)\dm_2^2+\dm v^m_t\dm_2^2\leq \dm v^m(t)\dm_\infty^2\dm\n v^m(t)\dm_2^2\,,\mbox{ for all }t>0\,.\ee  Assume that, for all $m\in\N$, $v^m$ solves a suitable system in such a way that the weak limit $v(t,x)$, weak limit with respect to the metric detected by the energy relation \rf{EW}, is a weak solution to the IBVP \rf{NS}. \par Then the weak solution $v(t, x)$ enjoys the following property: there exist $\theta\leq c\dm v_0\dm_2^4$ and a family of intervals $\{(\theta_\ell,T_\ell)\}_{\ell\in \A}$, $\A\subseteq \N$, such that $(0,\theta)-{\underset {m\in\A}\cup}(\theta_\ell,T_\ell)$ has zero Lebesgue's measure, and \be\label{ST}\ba{c}v\in C(\theta_\ell,T_\ell;J^{1,2}(\OO))\cap L^2(\theta_\ell,T_\ell;W^{2,2}(\OO))\mbox{\;  and\; }v_t\in L^2(\theta_\ell,T_\ell;J^2(\OO))\,,\VS v\in C(\theta,T;J^{1,2}(\OO))\cap L^2(\theta,T;W^{2,2}(\OO))\mbox{ \; and\; }v_t\in L^2(\theta,T;J^2(\OO))\,, \mbox{ for all }T>\theta.\ea\ee}\end{tho}
Of course, the Theorem\,\ref{CT} is interesting  if we are able to prove that there exists at least a sequence $\{v^m\}$ enjoying the hypotheses of the theorem. As we will point out, this is the case of the  Galerkin Hopf approximations constructed for the IBVP in a bounded domain $\OO$ and the 
sequence constructed by means of the Leray's approach in the setting of the  IBVP stated for \rf{NS}.\vskip0.15cm
{Now, assume that the definition of the Leray's weak solution $v$ is the one that considers the weak solution achieved as limit of the sequence $\{v^m\}$ of solutions to the mollified Navier-Stokes problem \rf{NS}. Then {\it a priori} for this weak solution we have two different {\it Th\'{e}or\{e}me de structure}. One is the classical Leray's structure theorem, another is the theorem proposed in this note. We  claim that the two {\it Th\'{e}or\`{e}me de structure} are not comparable  at most on a subset of $t>0$ with zero Lebesgue measure. Actually, the Leray {\it Th\'{e}or\`{e}me de structure} holds at most for all $t\in {\mathcal T}_L$ (see \rf{PLIM}) whose difference with $(0,\infty)$ is at most for a set of zero measure. Analogously, in Theorem\,\ref{CT} we consider the set $\mathscr T$ defined in \rf{SC-II}, that is different from $(0,\infty)$ at most for a subset of zero measure, and {\it a priori} $\mathscr I$ is not comparable with $\mathcal T$. Hence, the regularity of $v$, furnished by the two theorems, is not comparable at most for infants $t$ belonging to a subset of zero measure.
\vskip0.15cm We conclude the remark on the statement of Theorem\,\ref{CT} pointing out that,   in the forthcoming paper \cite{MP}, the authors employ the methods of proof of the present theorem to achieve a  {\it Th\'eor\�me de Structure} for a weak solution in the case of the dynamic interaction between a rigid body and a viscous fluid. 
\begin{rem}\label{R-I}{\rm Usually, the energy inequality \rf{ES}  is  deduced  by employing for $v^m(s,x)$   a   compactness theorem for Lebesgue $L^2$ space. For this goal, in the case of bounded domains $\OO$, one employs  the Rellich-Kondracev theorem of $W^{1,2}(\OO)\hookrightarrow L^2(\OO)$. In the case of an unbounded domain, one arrives at the compactness result  by using on $\OO\cap B_R$ the quoted theorem of compact immersion and by using an estimate of the uniform absolute continuity of the $L^2$-norm in the exterior of the ball $B_R$ for the sequence $\{v^m\}$, {\it e.g.} estimate \rf{CT-II}, that is achieved  by means of a suitable estimate, involving the sequence of pressure $\{\pi_{v^m}\}$ terms too. This last fact, in 3D, makes the difference between a Leray's weak solution and an Hopf's weak solution related to the IBVP in unbounded domains $\OO$.\par 
In our special approach to the structure theorem,  as it will be clear from our proof, in order to deduce, almost everywhere in $t>0$, the strong convergence  of $\{\n v^m\}$ to $\n v(t, x)$ in the $L^2$-norm, $v$ weak solution,  for the sequence $\{v^m\}$ we need the property of the uniform absolute continuity of the $L^2$-norm on the cylinder $(0,t)\times\{\OO-B_R\}$, that we achieve by employing \rf{CT-II}. Although with different aims, we claim that   between the Leray {\it Th\'eor\�me de Structure}  and our proof the unique common point   is  an estimate  of the kind \rf{CT-II}. Of course, in the case of the IBVP in bounded domains  $\OO$ estimate \rf{CT-II} is not employed.\par {We  point out that in our setting of a priori estimates related to the sequence $\{v^m\}$ in order to deduce partial regularity, although we provide a proof only for the {\it Th\'eor\�me de Structure}, one can consider also the extra assumptions of Prodi-Serrin type. These extra assumptions are usually employed on the weak solutions. {In analogy with the paper \cite{PM-PS},}      we can require the extra assumptions locally in time,  but, following the approach given in this note, with the difference that they are employed   on the elements of the sequence $\{v^m\}$. We give no details for this task.}}\end{rem}
\section{On the existence of a set of sequences of solutions enjoying properties \rf{CT-I}-\rf{DDN} of Theorem\,\ref{CT}}The aim of this section is to prove that    the  set of sequences $\{v^m\}$ enjoying the assumptions of Theorem\,\ref{CT} is not empty. 
\par {We first recall the following result.}
\begin{lemma}\label{LIPD}{Let $g\in J^{1,2}(\OO)\cap W^{2,2}(\OO)$. Then there exists a constant $c$, independent  of $g$, such that \be\label{LIPD-I}\dm g\dm_\infty\leq c\dm P\Delta g\dm_2^\frac12\dm \n g\dm_2^\frac12\,.\ee}\end{lemma}  \bp See Theorem\,1.1 in \cite{MRIV}.\ep
\par
We introduce the mollified Navier-Stokes system and the related IBVP:
\be\label{NSM}\ba{cc}v_t^m -\Delta v^m+\n\pi_{v^m}= -\mathbb J_m[v^m]\cdot\n v^m\,,&\n\cdot v^m=0\,,\mbox{ on }(0,T)\times \OO\,,\VS v^m=0\mbox{ on }(0,T)\times\po\,,& v^m=v^m_0\mbox{ on }\{0\}\times\OO\,,\ea\ee
where $\mathbb J_m[\,\cdot\,]\equiv J_{\frac1m}[\,\cdot\,]$ and $J_{\frac1m}[\,\cdot\,]$ is the Friedrichs (space) mollifier,  and $\{v_0^m\}\subset\mathscr C_0(\OO)$ converges to $v_0$ in the $L^2$-norm.\par Let $\OO$ be a smooth bounded domain and system  $\{a^p\}$ be an  basis in $J^2(\OO)\cap J^{1,2}(\OO)$, orthonormal in $J^2(\OO)$ where, for all $p\in\N$, we assume $a^p$ eigenvector of Stokes problem in $\OO$ with eigenvalue $\lambda_p$, see {\it e.g.} \cite{GPG}. We also introduce the Galerkin approximation in the sense of Prodi \cite{Pr} (see also \cite{Hy}). We set \be\label{GAP}u^n(t, x):=  \mbox{$\overset n{\underset{h=1}\sum}$}a^h(x)c^n_h(t)\ee being, for all $n\in\N$, the coefficients $\{c^n_h(t)\}$, for $h=1,\dots,n$, considered as solutions to the following system 
\be\label{NSG}\ba{c}(u^n_t,a^p)+(\n u^n,\n a^p)=-(u^n\cdot\n u^n,a^p)\,,\;p=1,\cdots, n\,,\VS u^n(0):=\mbox{$\overset n{\underset{p=1}\sum}$}(v_0,a^p)\,.\ea\ee     \par
The following lemmas hold:
\begin{lemma}\label{ENSM}{\sl For all $m$ there exists a regular solution $(v^m,\pi_{v^m})$ to problem \rf{NSM} defined for all $t>0$ and enjoying the properties:
\be\label{ENSM-I}v^m\in C(0,T;J^{1,2}(\OO))\cap L^2(0,T;W^{2,2}(\OO))\mbox{  and }v^m_t\in L^2(0,T;J^2(\OO))\,, \mbox{ for all }T>0\,,\ee and estimates \rf{EW}, \rf{DS} and \rf{DDN} hold{, with $\alpha=\frac23$.}}\end{lemma}
\bp The  proof  related to the existence and to the estimate \rf{EW} is well known, cf. {\it e.g.} \cite{GM,M-EE}. Let us consider  estimate \rf{DDN}. Actually, 
 for all $m\in\N$, after projecting   the equation \rf{NSM} by means of operator $P_2\equiv P$, {we evaluate the $L^2$-norm of both side terms, that is
\be\label{DN}\ba{ll}\sfrac d{d\tau}\dm \n v^m(\tau)\dm_2^2+\!\dm P\Delta v^m(\tau)\dm_2^2+\dm v^m_\tau\dm_2^2\hskip-0.25cm&=\!\dm P(\mathbb J_m[v^m(\tau)]\!\cdot\!\n v^m(\tau))\dm_2^2\leq\! \dm \mathbb J_m[v^m(\tau)]\!\cdot\!\n v^m(\tau)\dm_2^2\VSE\leq \dm v^m(\tau)\dm_\infty^2\dm\n v^m(\tau)\dm_2^2\,.\ea\ee
Finally, estimate \rf{DS} with $\alpha=\frac23$ is obtained in  \cite{CGM-I,DF}.  We believe that the proof given in \cite{CGM-I}  is  the simplest of  those quoted. We reproduce the proof. By virtue of estimate \rf{LIPD-I}, the right hand side of \rf{DN} is increased in the following way:
$$\sfrac d{d\tau}\dm \n v^m(\tau)\dm_2^2+\!\dm P\Delta v^m(\tau)\dm_2^2+\dm v^m_\tau\dm_2^2\leq c\dm P\Delta v^m(\tau)\dm_2\dm\n v^m(\tau)\dm_2^3\,.$$ Hence, trivially, we get
$$(1+\dm\n v^m(t)\dm_2^4)^{-1}\sfrac d{d\tau}\dm \n v^m(\tau)\dm_2^2+\sfrac12(1+\dm \n v^m(t)\dm_2^4)^{-1}\!\dm P\Delta v^m(\tau)\dm_2^2\leq c\dm\n v^m(\tau)\dm_2^2\,.$$ Integrating on $(0,T)$ and employing Holder's inverse inequality with exponent $-\frac12$ and $\frac13$, we arrive at
$$\sfrac12\Big[\intll0T\dm P\Delta v^m(t)\dm_2^\frac23dt\Big]^3\leq \Big[T+\intll0T\dm \n v^m(t)\dm_2^2dt\Big]^2\Big[\arctan(\dm \n v^m(0)\dm_2^2)+c \intll0T\dm \n v^m(t)\dm_2^2dt\Big]\,.$$
 The lemma is completely proved.}\ep \begin{lemma}\label{LDS}{\sl Let $\{(v^m,\pi_{v^m})\}$ be the sequence stated in Lemma\,\ref{ENSM}. Then for the pressure field we get
\be\label{RPL}\dm \pi_{v^m}(t)\dm_{L^2(|x|>R)}=o(1)\,\mbox{ for all }t>0\,.\ee}\end{lemma}
      \par The result of Lemma\,\ref{LDS} is also well known in literature. We limit ourselves to quote the one contained in \cite{CGM-I,GM,L}\,. 
\par For the Galerkin approximation we have the same results.
\par
\begin{lemma}\label{GL}{\sl For $u^n$ defined in \rf{GAP} it holds \be\label{GL-I}u^n\in C(0,T;J^{1,2}(\OO))\cap L^2(0,T;W^{2,2}(\OO))\mbox{  and }u^n_t\in L^2(0,T;J^2(\OO))\,, \mbox{ for all }T>0\,,\ee and estimates \rf{EW} and \rf{DS}.}\end{lemma}\bp 
Properties \rf{GL-I} and \rf{EW} are well known, cf {\it e.g.} \cite{Pr,Hy}.
We do not prove \rf{DS}. The estimate is analogous to the one given in \cite{CGM-I} for the sequence $\{v^m\}$ of solutions to problem \rf{NSM}. The unique difference is in the kind of approximation. \ep
\section{Proof of Theorem\,\ref{CT}}
\begin{lemma}\label{LSC}{\sl Let $\{(v^m,\pi_{v^m})\}$ ( resp. $\{u^{n}\}$) be the sequence enjoying \rf{EW}-\rf{DS} (of Lemma\,\ref{GL}). Denote by $v$ weak limit with respect the energy metric \rf{EW} ( resp.\,$u$ weak limit of $\{ u^n\}$). Then, almost everywhere in $t>0$, we get
$\displ\lim_m\dm \n v^m(t)-\n v(t)\dm_2=0$ ( resp.\,\,$\displ\lim_n\dm \n u^n(t)-\n u(t)\dm_2=0$)\,.}\end{lemma}
\bp In the case of the IBVP in bounded domain $\OO$,
the proof of the lemma is the same for both the sequences $\{v^m\}$ and $\{u^n\}$. Actually, applying the H\"older inequality, we get
$$\dm \n u^n(\tau)-\n u^p(\tau)\dm_2^2=-(P\Delta (u^n(\tau)-u^p(\tau)),u^n(\tau)-u^p(\tau))\leq \dm P\Delta (u^n(\tau)-u^p(\tau))\dm_2\dm u^n(\tau)-u^p(\tau)\dm_2 \,.$$ Integrating on $(0,T)$, we arrive at
\be\label{SPC}\intll0T\dm \n u^n(\tau)-\n u^p(\tau)\dm_2\leq \intll0T\dm P\Delta u^n(\tau)-P\Delta u^p(\tau)\dm_2^\frac12\dm u^n(\tau)-u^p(\tau)\dm_2^\frac12d\tau\,.\ee By virtue of \rf{DS} and the Friederich's lemma \cite{Ld}, that ensures the strong convergence of $\{u^n\}$ in $L^2(0,T;L^2(\OO))$, we achieve that  $\{\n u^n\}$ satisfies a Cauchy condition in $L^1(0,T;L^2(\OO))$. Therefore, the strong convergence in $L^1(0,T;L^2(\OO))$, that, finally, implies the thesis of the lemma. In the case of the IBVP in an unbounded domain $\OO$,  of course 
we limit ourselves to the case of the sequence $\{v^m\}$. We modify the estimate \rf{SPC}:
$$   \intll0T\!\dm \n v^m(\tau)-\n v^p(\tau)\dm_2\!\leq\! \intll0T\!\dm P\Delta v^m(\tau)-P\Delta v^p(\tau)\dm_2^\frac12\Big[\dm v^m(\tau)-v^p(\tau)\dm_{\null_{L^2(\OO\cap\{|x|<R\}}}^\frac12\!+2\dm v^m(\tau)\dm_{\null_{L^2(|x|>R)}}^\frac12 \Big]d\tau\,.$$ Now, the hypotheses of the theorem allow us to deduce again the Cauchy condition. Hence, we conclude the proof.
\ep \par{\it Proof of Theorem\,\ref{CT}.}\par
Taking energy estimate \rf{EW} into account, we are going to prove that for all $\eta>0$ there exists a $t^m\in[0,\eta^{-2}\dm v_0\dm_2^4)\equiv[0,\theta)$ such that
\be\label{UI-II}\dm v_0\dm_2\dm\n v^m(t^m)\dm_2\leq\eta\,,\mbox{ for all }m\in\N\,.\ee
We contradict estimate \rf{UI-II}. Hence, given $\eta>0$, assume that  there exists an index $m$ such that $\dm v_0\dm_2\dm\n v^m(t)\dm_2>\eta$ for  all $t\in[0,\eta^{-2}\dm v_0\dm_2^4)$\,.
Then the following holds:
$$\ba{ll}\dm v_0^m\dm_2^2\hskip-0.2cm&\displ=\eta^{-2}\dm v_0\dm_2^2\dm v_0^m\dm_2^2\dm v_0\dm_2^{-2}\eta^2< 2\eta^{-2}\dm v_0\dm_2^4\dm v_0\dm_2^{-2}\eta^2\VSE=2\dm v_0\dm_2^{-2}\eta^2\hskip-0.3cm\intll0{\eta^{-2}\dm v_0\dm_2^4}dt\;< 2\hskip-0.3cm\intll0{\eta^{-2}\dm v_0\dm_2^4}\hskip-0.3cm\dm\n v^m(t)\dm_2^2dt\,,\ea $$ in the last we increased thanks the absurd assumption. Of course, the last estimate contradicts the energy inequality \rf{EW}. Hence, estimate \rf{UI-II} is proved. Employing \rf{DDN} we also get
\be\label{DDN-I}\dm\n v^m(t)\dm_2^{-4}\sfrac d{dt}\dm\n v^m(t)\dm_2^2\leq c\dm \n v^m(t)\dm_2^2\,,\mbox{ for all }t>0\mbox{ and }m\in\N.\ee The last estimate, via an integration in time, furnishes
\be\label{DDN-II}\hskip-0.3cm\ba{ll}\dm \n v^m(t)\dm_2^2\hskip-0.3cm&\leq\! \dm\n v^m(t^m)\dm_2^2\Big[1\!-c\dm\n v^m(t^m)\dm_2^2\dm v^m(t^m)\dm_2^2\Big]^{-1}\VSE\leq\! \dm\n\hskip-0.025cm v^m(t^m)\dm_2^2\Big[1\!-\!c\dm\n\hskip-0.025cm v^m(t^m)\dm_2^2\dm v^m_0\dm_2^2\Big]^{-1}\hskip-0.2cm\leq \!\dm \n\hskip-0.025cm v^m(t^m)\dm_2^2\Big[1\!-\!2c\dm\n\hskip-0.025cm v^m(t^m)\dm_2^2\dm v_0\dm_2^2\Big]^{-1}\VSE\leq \eta^2\dm v_0\dm_2^{-2}\Big[1-2c\eta^2\Big]^{-1}\,,\mbox{ for all }t\geq \theta\geq t^m\mbox{ and for all }m\in\N\,,\ea \ee provided that $2c\eta^2<1$. Taking the last estimate in the hands, via an integration in time of estimate \rf{DDN}, we also deduce
\be\label{DDN-III}\intll\theta t\Big[\dm P\Delta v^m(\tau)\dm_2^2+\dm v^m_{\tau}(\tau)\dm_2^2\Big]d\tau\leq c(\dm v_0\dm_2,\eta)\,,\mbox{  for all }t\geq  \theta\geq t^m\mbox{ and for all }m\in\N\,.\ee
 Since the following holds: $$\mbox{\Large$\underset{n\in\N}\cap$}(t_m,\infty)\supseteq[\theta,\infty),$$ we realize that the sequence of regular solutions $\{(v^m,\pi_{v^m})\}$ admits
a common interval $(\theta,\infty)$ of  existence and of regularity like that detected by  the norms \rf{RS}.\par Now, we are going to consider the regularity for $t<\theta$.
 Employing the assumptions of the theorem,   Lemma\,\ref{LSC} ensures that    \be\label{SC-I}\mbox{a.e. in }t\in[0,\theta],\quad \dm \n v^m(t)\dm_2\to\dm \n v(t)\dm_2\,.\ee We set \be\label{SC-II}\mathscr T:=\{t:\mbox{the limit property 
\rf{SC-I} holds}\}\,.\ee We are going to prove that  for all $t\in \mathscr T$ there exists $m(t)$ such that, for all $m\geq m(t)$,  $(v^m,\pi_{v^m})$ admits a  maximal interval $(t,t+T(t,m))$ of existence and of regularity which contains strictly the following  $(t,t+c(\dm \n v(t)\dm_2^2+1)^{-2}]$\,.\par Actually,   
by virtue of estimate \rf{LIPD-I}, we estimate the right-hand side of \rf{DDN} as follows:
$$ \dm v^m(\tau)\dm_\infty^2\dm \n v^m(\tau)\dm_2^2\leq c\dm P\Delta v^m(\tau)\dm_2\dm \n v^m(\tau)\dm_2^3\leq \sfrac12\dm P\Delta v^m(\tau)\dm_2^2+c\dm\n v^m(\tau)\dm_2^6\,,$$ with a constant $c$ independent of $m$ and $\tau$. Hence, for estimate \rf{DDN} we arrive at
\be\label{ITG}\sfrac12\sfrac d{d\tau}\dm \n v^m(\tau)\dm_2^2+\sfrac12\dm P\Delta v^m(\tau)\dm_2^2+\dm v^m_\tau(\tau)\dm_2^2\leq 
c\dm\n v^m(\tau)\dm_2^6\,,\ee with a constant $c$ independent of $m$ and $\tau$. By virtue of \rf{SC-I}, for all $t\in\mathscr T$, there exists a $m(t)$
such that $\dm\n v^m(t)\dm_2\leq \dm \n v(t)\dm_2+1$, for all $m\geq m(t)$.
An integration of \rf{ITG} leads
to the estimates
\be\label{DDN-IV}\dm \n v^m(s)\dm_2^2+c\intll{t}s\Big[\dm P\Delta v^m(\tau)\dm_2^2+\dm v^m_\tau\dm_2^2\Big]d\tau<\infty\,,\mbox{ for all }s\in (t,t+T(t,m))\,,\ee  where $(t,t+T(t,m))$ is  a maximal interval  with property $(t,t+T(t,m))\supset (t, t+c(\dm\n v(t)\dm_2^2+1)^{-2})$, for all $m\geq m(t)$.
Considering $(t,t+T(t)):={\underset{m\geq m(t)}\cap}(t,t+T(m))$, for all $m\geq m(t)$, we realize a non-empty  maximal  interval of  existence and of regularity like that detected by  the norms \rf{RS}. \par By virtue of a classical argument one proves that the sequence $(v^m,\pi_{v^m})$ of the hypotheses of the theorem
admits a weak limit $v$ with respect to the metric of the energy \rf{EW}, with  $(v(t),\varphi)$ continuous function on $(0,\infty)$, for all $\varphi\in \mathscr C_0(\OO)$, and $v$  weak solution to problem \rf{NS} and $\displ \lim_{t\to0}\dm v(t)-v_0\dm_2=0\,.$
Of course, the previous properties continue to hold for the limit of any sequence extract from $\{v^m,\pi_{v^m}\}$. From estimates \rf{DDN-II} and \rf{DDN-III} we can extract a sequence with limit $w$ regular solution on $(\theta,\infty)$. On the other hand, due to the uniqueness of the weak limit in $L^2(\theta,\infty;J^{1,2}(\OO))$, we also  get that $v\equiv w$. This proves that the weak solution $v$, for $t>\theta$,  is regular in the sense of the Definition\,\ref{RS}.\par Now, let us consider\,\footnote{\,We essentially apply a procedure already employed in \cite{M-DD}. However, in \cite{M-DD} the author works on the weak solution and not on the approximating sequences of the weak solution.} $t_\alpha\in [0,\theta]\cap \mathscr T$. By virtue of \rf{DDN-IV} we can consider an extract that admits a limit $(w,\pi_w)$ regular solution on some maximal interval of the kind $I_\alpha:=(t_\alpha,t_\alpha+T(t_\alpha))$ containing strictly the interval $(t_\alpha,t_\alpha+c(\dm\n v(t_\alpha)\dm_2^2+1)^{-2}]$. By virtue of the uniqueness of the weak limit in  $L^2(0,\theta;J^{1,2}(\OO))$, we  get that $v\equiv w$ again. Hence, we have proved a local in time property of regularity for $v$. Let us consider $t_\beta\in [0,\theta]\cap \mathscr T-I_\alpha$. By the same previous argument lines,  we can state the existence of $I_\beta$ like maximal  interval of regularity for the weak solution $v$. Now, there are two possibilities:
\begin{itemize}\item[i.] $I_\alpha\cap I_\beta=\emptyset$;\item[ii.] $I_\alpha\cap I_\beta\ne\emptyset$, then $I_\beta\supset I_\alpha$.\end{itemize}
We justify the second item. If  $I_\alpha\cap I_\beta\ne\emptyset$, since $t_\beta\notin I_\alpha$, then $t_\beta<t_\alpha$. So that we have $t_\alpha\in I_\beta$. Hence $I_\alpha $ represents an extension of $I_\beta$, that is $t_\beta+T(t_\beta)=t_\alpha+T(t_\alpha)$. In the case of item i. We start again with a new $t_\gamma\in [0,\theta]\cap \mathscr T-\{I_\alpha\cup I_\beta\}$. In the case of item ii. we replace $I_\beta$ with $I_\alpha$ and start with $t_\gamma\in [0,\theta]\cap \mathscr T-I_\beta$. Iterating the procedure, we state a family of intervals $\{I_\alpha\}_{\alpha\in\mathscr A}$, with $\mathscr A$ set of indexes. We have that $[0,\theta]\cap \mathscr T\subseteq {\underset{\alpha\in \mathscr A}\cup}I_\alpha$, with ${\overset\circ{I_\alpha}}\ne\emptyset$ and ${\overset\circ{I_\alpha}}\cap{\overset\circ{ I_\beta}}=\emptyset$ for $\alpha\ne\beta$. This last property ensures that the set of indexes $\mathscr A$ has at most the cardinality of $\N$. The theorem is completely proved.

\end{document}